\documentclass[11pt]{article}
\usepackage{graphicx,subfigure}
\usepackage{epsfig}
\usepackage{amssymb}
\usepackage{mathrsfs}
\usepackage{color}
\usepackage[dvipdfm]{hyperref}

\newcommand{\E}{{\cal E}}

\newtheorem{theorem}{Theorem}[section]

\def\whitebox{{\hbox{\hskip 1pt
 \vrule height 6pt depth 1.5pt
 \lower 1.5pt\vbox to 7.5pt{\hrule width
    3.2pt\vfill\hrule width 3.2pt}%
 \vrule height 6pt depth 1.5pt
 \hskip 1pt } }}
\def\qed{\ifhmode\allowbreak\else\nobreak\fi\hfill\quad\nobreak
     \whitebox\medbreak}

\newcommand{\ignore}[1]{}

\begin {document}

\baselineskip 16pt
\title{Asymptotic incidence energy and Laplacian-energy-like\\ invariant of
the Union Jack lattice}

 \author{\small   Jia\textrm{-}Bao \ Liu$^{a,b}$, \ \  Xiang\textrm{-}Feng \ Pan$^{a,}$\thanks{Corresponding author. Tel:+86-551-63861313. \E-mail:liujiabaoad@163.com (J.Liu), xfpan@ahu.edu.cn(X.Pan).}\\
 \small  $^{a}$ School of Mathematical Sciences, Anhui  University, Hefei 230601, P. R. China\\
 \small  $^{b}$ Department of Public Courses, Anhui Xinhua
 University, Hefei 230088, P. R. China\\}

\date{}
\maketitle
\begin{abstract}
  The incidence energy
$\mathscr{IE}(G)$ of a graph $G$, defined as the sum of the
singular values of the incidence matrix of a graph $G$, is a much
studied quantity with well known applications in chemical physics.
The Laplacian-energy-like invariant of $G$ is defined as the sum
of square roots of the Laplacian eigenvalues.
 In this paper, we obtain the closed-form formulae
expressing the incidence energy and the Laplacian-energy-like
invariant of the Union Jack lattice. Moreover, the explicit
asymptotic values of these quantities
 are calculated by utilizing the applications of analysis approach with
the help of calculational software.

\medskip
\noindent {\bf Keywords}: Solvable lattice models;\ \ Random
graphs;\ \ Networks;  \  \  Incidence energy; \ \
Laplacian-energy-like invariant
\end{abstract}

\section{ Introduction}

 Lattices have several attractive features that make them
interesting candidates for use in matter physics.  The quantum
spin model with frustration and the Ising model of the Union Jack
lattice have been exploited extensively by
physicists~\cite{LY2015}.

The problem of calculation of some physical and chemical indices
(such as the energy, the incidence energy and the
Laplacian-energy-like invariant) on the lattices has been
extensively studied and it became a popular topic of research in
mathematical chemistry and mathematics. The energy of a graph $G$
arising in chemical physics, is defined as the sum of the absolute
values of the eigenvalues of $G$. The energy of many lattices were
considered by physicists~\cite{LY2015,LY,YLZ,Liu2014}.
 A general problem of interest in
physics, chemistry and mathematics is the calculations of various
energies of lattices~\cite{LY2015,YZ,YLZ,WY,ZZ}. The closed-form
formulae expressing the incidence energy of the hexagonal lattice,
triangular lattice, and $3^3.4^2$ lattice are investigated
in~\cite{Liu2015}. The Laplacian-energy-like invariant of
hexagonal, triangular, and $3^3.4^2$ lattices with three boundary
conditions are reported in~\cite{LiuP2015}.

  The authors
of~\cite{LY2015} investigated the formulae of the number of
spanning trees, the energy, and the Kirchhoff index of the Union
Jack lattice with toroidal boundary condition. Although the
incidence energy and the Laplacian-energy-like invariant of some
lattices are reported by mathematics and physics journals (see for
example Refs. \cite{Liu2015,LiuP2015}). But as far as we know, no
one has considered the incidence energy and the
Laplacian-energy-like invariant of the Union Jack lattice. In this
paper, we obtain the solution of these problems.

The rest of the paper is organized as follows. We introduce the
preliminaries and the definitions of lattices in Sections 2 and 3,
respectively. In section 4, we deduce the signless Laplacian
eigenvalues and the incidence energy of $UJL(n,m)$. In section 5,
we investigate the Laplacian-energy-like invariant of the Union
Jack lattice.

\section{Preliminaries}
 At the beginning of this section, we first
introduce some notations, which will be used in following
discussion.

 Let $G = (V(G),E(G))$ be a graph with vertex set $V(G)
= \{v_1, v_2, \dots , v_n\}$ and edge set $E(G) = \{e_1, e_2,
\dots , e_m\}$. Denote by $\mid V(G)\mid$ and $\mid E(G)\mid$ the
numbers of vertices and edges, respectively. The adjacency matrix
of graph $G$ is an $n \times n$ $(0,1)$-matrix $A(G) = (a_{ij})$,
where $a_{ij} = 1$ if and only if $(v_i, v_j)$ is an edge of $G$
and $a_{ij} = 0$ otherwise. Let $\lambda_1(G)\geq \lambda_2(G)\geq
\dots \geq \lambda_n(G)$ be the eigenvalues of the adjacency
matrix $A(G).$
$$ Spec_A(G) = \Big\{\lambda_1(G), \lambda_2(G), \dots
, \lambda_n(G)\Big\}$$
 is also called the spectrum of
$G$~\cite{CDS}.

The degree of a vertex $v$, denoted by $d_G(v)$, is the number of
edges incident to $v$ in a graph $G$. Let $D(G)$ be the diagonal
matrix of vertex degrees of $G$. The Laplacian matrix of $G$ is
$L(G)=D(G)-A(G)$ and the Laplacian spectrum of $G$ is denoted by
$$Spec_L(G) = \Big\{\mu_1(G), \mu_2(G), \dots , \mu_n(G)\Big\},$$
where $\mu_1(G), \mu_2(G), \dots ,
\mu_n(G)$ are the eigenvalues of the Laplacian matrix $L(G).$
 The signless Laplacian matrix is $Q(G)=D(G)+A(G)$ and the signless
Laplacian spectrum of $G$ is denoted by
$$ Spec_Q(G) = \Big\{q_1(G), q_2(G), \dots , q_n(G)\Big\},$$
 where $q_1(G), q_2(G), \dots , q_n(G)$
are the eigenvalues of the signless Laplacian matrix $Q(G).$
 It is well known that $L(G)$ and $Q(G)$
are symmetric and positive semi-definite, and the spectra of
$L(G)$ and $Q(G)$ coincide if and only if the graph G is
bipartite~\cite{CRS,DC}.

 Let $I(G)$ be the (vertex-edge) incidence matrix of the
graph $G$. The $(i,j)$-entry of $I(G)$ is 1 if $v_i$ is incident
with $e_j$ and 0 otherwise. (In what follows, the unit matrix of
order $n$ will be denoted by $E_n$ to avoid confusion with the
incidence matrix.) A well known fact is the identity~\cite{CDS}:
$$I(G) I(G)^t = A(G)+ D(G)= Q(G).$$

We recall Cartesian product of graphs and Kronecker product, which
will be used in the proof of our result.

 Given graphs $G$ and $H$ with vertex sets $U$ and $V$, the Cartesian product
 $G\Box H$ of graphs $G$ and $H$ is a graph such that
the vertex set of $G\Box H$ is the Cartesian product $U\Box V$;
and any two vertices $(u,u')$ and $(v,v')$ are adjacent in $G\Box
H$ if and only if either $u = v$ and $u'$ is adjacent $v'$ in $H$,
or $u' = v'$ and $u$ is adjacent $v$ in $G$ {\rm \cite{CZ}.

The Kronecker product $A\otimes B$ of two matrices $A =
(a_{ij})_{m\times n}$ and $B = (b_{ij})_{p\times q}$ is the
$mp\times nq $ matrix obtained from $A$ by replacing each element
$a_{ij} $ by $a_{ij}B $. If $A, B, C$ and $D$ are matrices of such
size that one can form the matrix products $AC$ and $BD$, then
$(A\otimes B)(C\otimes D)=AC\otimes BD.$ It follows that $A\otimes
B$ is invertible if and only if $A$ and $B$ are invertible, in
which case the inverse is given by $(A\otimes
B)^{-1}=A^{-1}\otimes B^{-1}.$ Note also that $(A\otimes
B)^{T}=A^{T}\otimes B^{T}$. Moreover, if the matrices $A$ and $B$
are of order $n\times n$ and $p\times p$, respectively, then
$det(A\otimes B) = (det~A)^p(det~B)^n.$ The readers are referred
to~\cite{Horn1991} for other properties of the Kronecker product
not mentioned here.

\section{ The $4.8.8$ lattice and the Union Jack lattice}

Our notations for the 4.8.8 lattice and the Union Jack lattice
follow~\cite{LY2015}. The 4.8.8 lattice with toroidal boundary
conditions, denoted by $G^t(n,m)$, is illustrated in Figure 1 (a).
The left and right (resp. the lower and upper) boundaries of the
picture are identified such that
  all $a_i$s, $a_i^*$s, $b_i$s, $b_i^*$s, are some vertices on
the left, right, lower and upper boundaries, respectively, and
$(a_1, a_1^*), (a_2, a_2^*),\dots , (a_m, a_m^*)$ and $(b_1,
b_1^*), (b_2, b_2^*), \dots , (b_n, b_n^*)$ are edges in
$G^t(n,m)$. Obviously, the 4.8.8 lattice $G^t(n,m)$ is composed of
$mn$ quadrangles.

\begin{figure}[ht]
\centering
  \includegraphics[width=\textwidth]{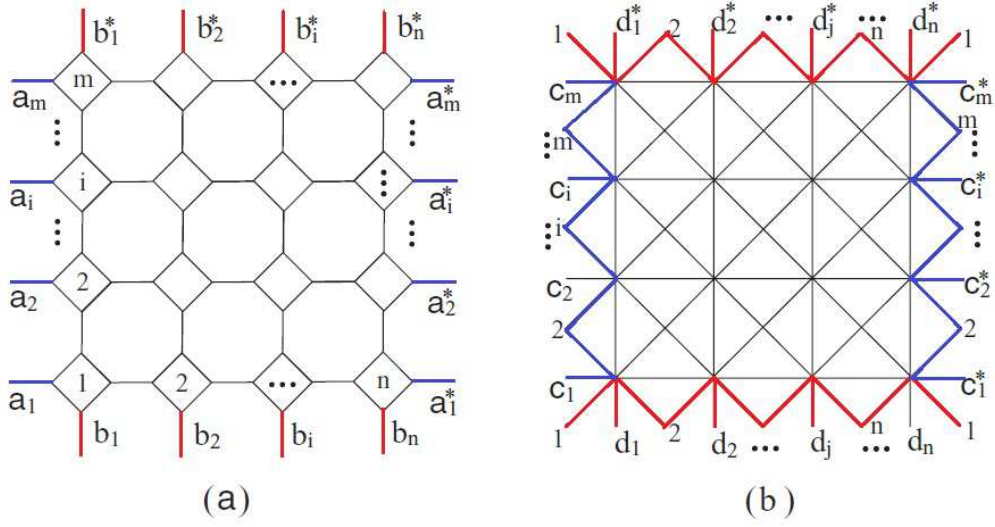}
\caption{(a) The 4.8.8 lattice $G^t(n,m)$ with toroidal boundary
condition.\ (b) The Union Jack lattice $UJL(n,m)$ with toroidal
boundary condition, i.e., the dual graph of the 4.8.8 lattice
$G^t(n,m)$.} \vspace{+1em}
\end{figure}

 The Union Jack lattice with toroidal boundary
condition, denoted by $UJL(n,m)$, is the dual lattice of the 4.8.8
lattice with toroidal boundary condition. Figure 1 (b) is the
Union Jack lattice $UJL(n,m)$ corresponding to the 4.8.8 lattice
illustrated in Figure 1 (a). Similarly, the left and right (resp.
the lower and upper) boundaries of the picture are identified such
that all $c_i$s, $c_i^*$s, $d_i$s, $d_i^*$s, are some vertices on
the left, right, lower and upper boundaries, respectively, and
$(c_1, c_1^*), (c_2, c_2^*),\dots , (c_m, c_m^*)$ and $(d_1,
d_1^*), (d_2, d_2^*), \dots , (d_n, d_n^*)$ are edges in
$UJL(n,m)$.

The Union Jack lattice with toroidal boundary condition $UJL(n,m)$
also can be obtained from an $n\times m$ square lattice with
toroidal boundary condition by inserting a new vertex $v_f$ to
each face $f$ and adding four edges $(v_f , u_i(f)), i=1,2,\dots ,
4$, where $u_i(f)$ are four vertices on the boundary of $f$. Let
$G'$ be the graph obtained from $UJL(n,m)$ by deleting the
vertices of degree 4, i.e., $G'$ is an $n\times m$ square lattice
with toroidal boundary condition, which is the Cartesian product
$C_n \Box C_m$ of two cycles $C_n$ and $C_m$. It is easy to see
that $UJL(n,m)$ has $2mn$ vertices.

\section{The signless Laplacian eigenvalues and the incidence energy of $UJL(n,m)$}

\subsection{The signless Laplacian eigenvalues of $UJL(n,m)$.}

In this subsection,
 we will deduce the signless Laplacian
eigenvalues of $UJL(n,m)$  by utilizing the techniques
in~\cite{LY2015}. For the convenience to following description, we
set
$$\alpha_i=\frac{2\pi i}{n}, \beta_j=\frac{2\pi j}{m},
i=0,1,\dots,n-1;j=0,1,\dots,m-1.$$

\begin{theorem}\label{4-1}
Let $UJL(n,m)$ be a Union Jack lattice with toroidal boundary
condition. Then the signless Laplacian eigenvalues of $UJL(n,m)$
are
$$ (6+ \cos\alpha_i+\cos\beta_j)\pm \sqrt{(6+ \cos\alpha_i+\cos\beta_j)^2-4(7 + \cos\alpha_i+\cos\beta_j-cos\alpha_i\cos\beta_j)},$$
where $0\leq i\leq n-1, 0\leq j\leq m-1.$
\end{theorem}

\noindent{\bf Proof.}
 With a suitable labelling of vertices of $UJL(n,m)$, the adjacency
matrix of $UJL(n,m)$ is
$$A(UJL(n,m)) = \left(%
\begin{array}{cc}
  A(C_n \Box C_m) & M \\
  M^T & 0 \\
\end{array}%
\right),$$ where $A(C_n \Box C_m)$ is the adjacency matrix of $C_n
\Box C_m$, $M$ is the matrix induced by the adjacency relation
between $V (C_n \Box C_m)$ and $V (F)$, here $F$ is the face set
of $C_n \Box C_m$, $M^T$ is the transpose of $M$. Moreover, $M$ is
a block matrix (with m times m entries) which has the following
form:
$M=\left(%
\begin{array}{cccccc}
  R & 0 & 0 &               \dots & 0 & R \\
  R & R & 0 & \dots & 0 & 0 \\
  0 & R & R &             \dots & 0 & 0 \\
  \vdots & \vdots & \ddots & \ddots & \ddots & \vdots \\
  0 & 0 & 0 & \dots & R & 0 \\
  0 & 0 & 0 & \dots & R & R \\
\end{array}%
\right)_{m\times m},$
where
$R=\left(%
\begin{array}{cccccc}
  1 & 0 & 0 &               \dots & 0 & 1 \\
  1 & 1 & 0 & \dots & 0 & 0 \\
  0 & 1 & 1 &             \dots & 0 & 0 \\
  \vdots & \vdots & \ddots & \ddots & \ddots & \vdots \\
  0 & 0 & 0 & \dots & 1 & 0 \\
  0 & 0 & 0 & \dots & 1 & 1 \\
\end{array}%
\right)_{n\times n}.$

On one hand, the degrees of $UJL(n,m)$ are 4 or 8. With a suitable
labelling of vertices of $UJL(n,m)$, the diagonal matrix
$D(UJL(n,m))$ is
\begin{equation}\label{}
D(UJL(n,m))=diag\textbf{\{}\underbrace{4,4,\dots,4}_{mn},\underbrace{8,8,\dots,8}_{mn}\textbf{\}}.
\end{equation}

On the other hand, the adjacency matrix $A(C_n \Box C_m)$ of $C_n
\Box C_m$ has the following form by a suitable labelling of
vertices of $C_n \Box C_m$:
\begin{eqnarray} \nonumber
 A\Big(C_n \Box C_m\Big)&=&\left(%
\begin{array}{cccccc}
  A(C_{n}) & E_{n} & 0 &               \dots & 0 & E_{n} \\
  E_{n} & A(C_{n}) & E_{n} & \dots & 0 & 0 \\
  0 & E_{n} & A(C_{n}) &             \dots & 0 & 0 \\
  \vdots & \vdots & \ddots & \ddots & \ddots & \vdots \\
  0 & 0 & 0 & \dots & A(C_{n}) & E_{n} \\
  E_{n} & 0 & 0 & \dots & E_{n} & A(C_{n}) \\
\end{array}%
\right)_{n\times n}\\
 &=& E_m \otimes A(C_n) + A(C_m) \otimes E_n, \\\nonumber
\end{eqnarray}
where $E_n$ is the identity matrix of order $n$.

Based on equations (1) and (2),
 the signless Laplacian matrix $Q(UJL(n,m))$ has the following form:

\begin{eqnarray} \nonumber
Q(UJL(n,m)) &=& D(UJL(n,m))+ A(UJL(n,m))\\
\nonumber  &=&\left(%
\begin{array}{cc}
  8E_{mn}+A(C_n \Box C_m) & M \\
  M^T & 4E_{mn} \\
\end{array}%
\right), \\\nonumber
\end{eqnarray}

Hence the signless Laplacian characteristic polynomial of
$UJL(n,m)$ is

\begin{eqnarray} \nonumber
\psi(UJL(n,m),x) &=& det~(xE_{2mn}- Q(UJL(n,m))\\
\nonumber &=& det~ \left(%
\begin{array}{cc}
  (x-8)E_{mn}-A(C_n \Box C_m) & -M \\
  -M^T & (x-4)E_{mn} \\
\end{array}%
\right) \\
 &=& det~\Big((x-4)(x-8)E_{mn}- (x-4)A(C_n \Box C_m)-MM^T\Big).\\
\nonumber
\end{eqnarray}

Considering that

\begin{eqnarray} \nonumber
MM^T&=&\left(%
\begin{array}{cccccc}
  2RR^T & RR^T & 0 &               \dots & 0 & RR^T\\
  RR^T & 2RR^T & RR^T & \dots & 0 & 0 \\
  0 & RR^T & 2RR^T &             \dots & 0 & 0 \\
  \vdots & \vdots & \ddots & \ddots & \ddots & \vdots \\
  0 & 0 & 0 & \dots & 2RR^T & RR^T \\
  RR^T & 0 & 0 & \dots & RR^T & 2RR^T \\
\end{array}%
\right)\\\nonumber
 &=&B_m\otimes B_n.
\end{eqnarray}
 where
$$B_n=\left(%
\begin{array}{cccccc}
  2 & 1 & 0 &               \dots & 0 & 1 \\
  1 & 2 & 1 & \dots & 0 & 0 \\
  0 & 1 & 2 &             \dots & 0 & 0 \\
  \vdots & \vdots & \ddots & \ddots & \ddots & \vdots \\
  0 & 0 & 0 & \dots & 2 & 1 \\
  1 & 0 & 0 & \dots & 1 & 2 \\
\end{array}%
\right)_{n\times n} =2E_n + A(C_n).$$

Consequently,
\begin{equation}\label{}
    MM^T = \Big(2E_m + A(C_m)\Big) \otimes \Big(2E_n + A(C_n)\Big).
\end{equation}

 According to
equations (2), (3) and (4), we have

\begin{eqnarray} \nonumber
\psi(UJL(n,m),x) &=&  det~\Big((x-4)(x-8)E_{mn}- (x-4)A(C_n \Box C_m)-MM^T\Big)\\
\nonumber &=& det~\Bigg\{(x-4)(x-8)E_{mn}- (x-4) \Big[E_m \otimes
A(C_n) + A(C_m)\otimes E_n \Big]\\ \nonumber
 &&- \Big(2E_m + A(C_m)\Big) \otimes \Big(2E_n + A(C_n)\Big)\Bigg\} \\
\nonumber &=& det~\Bigg\{(x^2-12x+28)E_{mn}- (x-2) \Big[E_m
\otimes A(C_n) + A(C_m)\otimes E_n \Big]\\ \nonumber
 &&- A(C_m) \otimes  A(C_n)\Bigg\} .\\
\nonumber
\end{eqnarray}

Let
\begin{equation}\label{}
Z=(x^2-12x+28)E_{mn}-(x-2)\big[E_m\otimes A(C_n)+A(C_m)\otimes
E_n\big]-A(C_m)\otimes A(C_n).
\end{equation}

Note that the eigenvalues of $A(C_n)$ are $2 \cos\alpha_i, 0\leq
i\leq n-1.$

 Hence, there exist two invertible matrices $P$ and $Q$
such that
$$P^{-1}A(C_n)P=diag\Big(2,2 \cos\alpha_i,\dots, 2 \cos\alpha_{n-1}\Big):=W_n,$$
and
$$Q^{-1}A(C_m)Q=diag\Big(2,2 \cos\alpha_i,\dots, 2 \cos\alpha_{m-1}\Big):=W_m.$$

 In fact, $P$ and $Q$ are invertible matrices, then $Q
\otimes P$ is an invertible matrix.

By the equation (5), one can obtain

\begin{eqnarray} \nonumber
(Q\otimes P)^{-1}Z(Q\otimes P) &=& (Q\otimes P)^{-1}
\Bigg\{(x^2-12x+28)E_{mn}- (x-2) \Big[E_m \otimes A(C_n)
\\\nonumber &&+
A(C_m)\otimes E_n \Big]- A(C_m) \otimes  A(C_n)\Bigg\} (Q\otimes P)\\
\nonumber  &=&  (x^2-12x+28)E_{mn}-(x-2)\Big[E_m\otimes
W_n+W_m\otimes E_n\Big]-W_m\otimes W_n. \\\nonumber
\end{eqnarray}
It is not difficult to see that
$(x^2-12x+28)E_{mn}-(x-2)\Big[E_m\otimes W_n+W_m\otimes
E_n\Big]-W_m\otimes W_n$
 is a diagonal matrix whose diagonal entries are
$(x^2-12x+28)-(x-2)( 2 \cos\alpha_i+2
\cos\beta_j)-4\cos\alpha_i\cos\beta_j. $

Actually, the zeros of $(x^2-12x+28)-(x-2)( 2 \cos\alpha_i+2
\cos\beta_j)-4\cos\alpha_i\cos\beta_j=0 $ are
$$ (6+ \cos\alpha_i+\cos\beta_j)\pm \sqrt{(6+ \cos\alpha_i+\cos\beta_j)^2-4(7 + \cos\alpha_i+\cos\beta_j-cos\alpha_i\cos\beta_j  )} .$$

Then the proof of this theorem is complete.  \qed

\subsection{The incidence energy of $UJL(n,m)$}

Gutman introduced the concept of energy $\mathscr{E}(G)$~{\rm
\cite{IG} for a simple graph $G$, which is defined as
$$\mathscr{E}(G)= \sum_{i=1}^{n}|\lambda_i(G)|.$$ As an analogue
to $\mathscr{E}(G)$, the incidence energy $\mathscr{IE}(G)$, is a
novel topological index.
 Inspired by Nikiforov's
idea~\cite{VN}, in 2009 Jooyandeh et al. \cite{JKM} introduced the
concept of incidence energy $\mathscr{IE}(G)$ of a graph $G$,
defining it as the sum of the singular values of the incidence
matrix $I(G),$ i.e.,
$$\mathscr{IE}(G)= \sum_{i=1}^{n}\sigma_i,$$
 where $\sigma_1,\sigma_2,\dots,\sigma_n$ are the singular values of the
incidence matrix $I(G)$.

Since the identity~\cite{CDS}: $I(G)I(G)^t=A(G)+D(G)=Q(G),$
 its immediate consequence is $$\sigma_i =\sqrt{ q_i}.$$
 Therefore $$\mathscr{IE}(G)= \sum_{i=1}^{n}\sqrt{q_i}.$$

 For more work on $\mathscr{IE}(G)$, the
readers are referred to~{\rm \cite{GKM,GKMZ,DG} and recent
articles ~\cite{SB,SBI,OR,TZ,ZJ,ZJZ,LC2014}.

 Next, we will explore the incidence energy
$\mathscr{IE}\Big(UJL(n,m)\Big)$ of $UJL(n,m).$

\begin{theorem}\label{4-2}
 Let $UJL(n,m)$ be a Union Jack lattice with toroidal boundary
condition. Then

1. The incidence energy $\mathscr{IE}\Big(UJL(n,m)\Big)$ of
$UJL(n,m)$ can be expressed by

\begin{eqnarray}
\nonumber  &&\mathscr{IE}\Big(UJL(n,m)\Big)\\
\nonumber &=& \sum_{i=0}^{n-1}\sum_{j=0}^{m-1}
\sqrt{(6+ \cos\alpha_i+\cos\beta_j)    + \sqrt{(6+ \cos\alpha_i+\cos\beta_j)^2-4(7 + \cos\alpha_i+\cos\beta_j-cos\alpha_i\cos\beta_j  )}}\\
\nonumber     &&+  \sum_{i=0}^{n-1}\sum_{j=0}^{m-1} \sqrt{(6+
\cos\alpha_i+\cos\beta_j)    -\sqrt{(6+
\cos\alpha_i+\cos\beta_j)^2-4(7 +
\cos\alpha_i+\cos\beta_j-cos\alpha_i\cos\beta_j  )}},\\ \nonumber
\end{eqnarray}
where $\alpha_i=\frac{2\pi i}{n}, \beta_j=\frac{2\pi j}{m},
i=0,1,\dots,n-1;j=0,1,\dots,m-1.$

2.  As $m,n\to \infty$, $\mathscr{IE}\Big(UJL(n,m)\Big)\approx
9.4770mn.$
\end{theorem}

\noindent{\bf Proof.} By Theorem 4.1 and the definition of the
incidence energy $\mathscr{IE}(G)$, one can obtain the incidence
energy $\mathscr{IE}\Big(UJL(n,m)\Big)$ of $UJL(n,m)$ is

\begin{eqnarray}
\nonumber  &&\mathscr{IE}\Big(UJL(n,m)\Big)\\
\nonumber &=& \sum_{i=0}^{n-1}\sum_{j=0}^{m-1}
\sqrt{(6+ \cos\alpha_i+\cos\beta_j)    + \sqrt{(6+ \cos\alpha_i+\cos\beta_j)^2-4(7 + \cos\alpha_i+\cos\beta_j-cos\alpha_i\cos\beta_j  )}}\\
\nonumber     &&+  \sum_{i=0}^{n-1}\sum_{j=0}^{m-1} \sqrt{(6+
\cos\alpha_i+\cos\beta_j)    -\sqrt{(6+
\cos\alpha_i+\cos\beta_j)^2-4(7 +
\cos\alpha_i+\cos\beta_j-cos\alpha_i\cos\beta_j  )}}.\\ \nonumber
\end{eqnarray}

Therefore the statement 1 of Theorem 4.2  is immediate.

In what follows, we will calculate the asymptotic value of the
incidence energy $\mathscr{IE}\Big(UJL(n,m)\Big)$. For the sake of
simplicity , we set
$$A=(6+ \cos\alpha_i+\cos\beta_j) +\sqrt{(6+ \cos\alpha_i+\cos\beta_j)^2-4(7 +\cos\alpha_i+\cos\beta_j-cos\alpha_i\cos\beta_j  )},$$
$$B=(6+ \cos\alpha_i+\cos\beta_j) -\sqrt{(6+ \cos\alpha_i+\cos\beta_j)^2-4(7 +\cos\alpha_i+\cos\beta_j-cos\alpha_i\cos\beta_j  )}.$$

Considering that $m,n$ approach infinity, we have
\begin{eqnarray}
\nonumber && \lim_{m\to \infty}\lim_{n\to \infty}\frac{
\mathscr{IE}\Big(UJL(n,m)\Big)}{2mn} \\
&=&\frac{1}{8\pi^2} \int_{0}^{2\pi}\int_{0}^{2\pi} \sqrt{A} \cdot
\ dx dy +\frac{1}{8\pi^2}\int_{0}^{2\pi}\int_{0}^{2\pi}
\sqrt{B} \cdot \ dx dy \\
   &\approx & 4.7385.
\end{eqnarray}

Consequently, according to the equality (7), we can get the
asymptotic incidence energy

$\mathscr{IE}\Big(UJL(n,m)\Big)\approx 9.4770mn,$ as $m,n\to
\infty$ . \qed

\vspace{5pt} \noindent {\bf Remark 4.3}\ The numerical integration
value in equality (6) is calculated with MATLAB software
calculation, i.e.,
$$\frac{1}{8\pi^2}\int_{0}^{2\pi}\int_{0}^{2\pi} \sqrt{A}\cdot\ dx
dy\approx2.9040,~~~ \frac{1}{8\pi^2}\int_{0}^{2\pi}\int_{0}^{2\pi}
\sqrt{B}\cdot\ dx dy\approx1.8345. $$

\section{The Laplacian-energy-like invariant of the Union Jack lattice}

 The Laplacian-energy-like invariant of a graph $G$, $\mathscr{LEL}(G)$ for short, is
defined as $$\mathscr{LEL}(G)=\sum_{i=1}^{n-1}\sqrt{\mu_i},$$
which is a novel topological index.
 The
concept of $\mathscr{LEL}(G)$ was first introduced by J. Liu and
B. Liu ({\rm \cite{LL}, 2008), where it showed that
$\mathscr{LEL}(G)$ has similar features as the graph energy
$\mathscr{E}(G)$~\cite{IZB}.

We recall the Laplacian eigenvalues of $UJL(n,m)$.  The authors of
\rm \cite{LY2015} proved the following theorem, whose proof thus
is omitted.

\begin{theorem}\label{5-1} (\cite{LY2015})
Let $UJL(n,m)$ be a Union Jack lattice with toroidal boundary
condition. Then the Laplacian eigenvalues of $UJL(n,m)$ are
$$ (6- \cos\alpha_i-\cos\beta_j)\pm \sqrt{(6- \cos\alpha_i-\cos\beta_j)^2-4(7 -3\cos\alpha_i-3\cos\beta_j-cos\alpha_i\cos\beta_j)},$$
where $0\leq i\leq n-1, 0\leq j\leq m-1.$
\end{theorem}

 The following theorem expresses the Laplacian-energy-like
invariant of the Union Jack lattice.

\begin{theorem}\label{5-2}
 Let $UJL(n,m)$ be a Union Jack lattice with toroidal boundary
condition. Then

1. The Laplacian-energy-like invariant
$\mathscr{LEL}\Big(UJL(n,m)\Big)$ of $UJL(n,m)$ can be expressed
by

\begin{eqnarray}
\nonumber  &&\mathscr{LEL}\Big(UJL(n,m)\Big)\\
\nonumber &=& \sum_{i=0}^{n-1}\sum_{j=0}^{m-1}
\sqrt{(6-\cos\alpha_i-\cos\beta_j)    + \sqrt{(6-\cos\alpha_i-\cos\beta_j)^2-4(7 -3\cos\alpha_i-3\cos\beta_j-cos\alpha_i\cos\beta_j  )}}\\
\nonumber     &&+  \sum_{i=0}^{n-1}\sum_{j=0}^{m-1}
\sqrt{(6-\cos\alpha_i-\cos\beta_j) -\sqrt{(6-
\cos\alpha_i-\cos\beta_j)^2-4(7 -3
\cos\alpha_i-3\cos\beta_j-cos\alpha_i\cos\beta_j  )}},\\ \nonumber
\end{eqnarray}
where $\alpha_i=\frac{2\pi i}{n}, \beta_j=\frac{2\pi j}{m},
i=0,1,\dots,n-1;j=0,1,\dots,m-1.$

2.  As $m,n\to \infty$, $\mathscr{LEL}\Big(UJL(n,m)\Big)\approx
9.3682 mn.$
\end{theorem}

\noindent{\bf Proof.} By Theorem 5.1 and the definition of the
Laplacian-energy-like invariant $\mathscr{LEL}(G)$, one can obtain
the $\mathscr{LEL}\Big(UJL(n,m)\Big)$ of $UJL(n,m)$ is

\begin{eqnarray}
\nonumber  &&\mathscr{LEL}\Big(UJL(n,m)\Big)\\
\nonumber &=& \sum_{i=0}^{n-1}\sum_{j=0}^{m-1}
\sqrt{(6-\cos\alpha_i-\cos\beta_j)    + \sqrt{(6-\cos\alpha_i-\cos\beta_j)^2-4(7 -3 \cos\alpha_i-3\cos\beta_j-cos\alpha_i\cos\beta_j  )}}\\
\nonumber     &&+  \sum_{i=0}^{n-1}\sum_{j=0}^{m-1}
\sqrt{(6-\cos\alpha_i-\cos\beta_j)
-\sqrt{(6-\cos\alpha_i-\cos\beta_j)^2-4(7-3\cos\alpha_i-3\cos\beta_j-cos\alpha_i\cos\beta_j
)}}.\\ \nonumber
\end{eqnarray}

Hence the statement 1 of Theorem 5.2 holds.

Similarly, we will formulate the asymptotic value of the
 $\mathscr{LEL}\Big(UJL(n,m)\Big)$.
Let
$$C=(6-\cos\alpha_i-\cos\beta_j) +\sqrt{(6-\cos\alpha_i-\cos\beta_j)^2-4(7-3\cos\alpha_i-3\cos\beta_j-cos\alpha_i\cos\beta_j  )},$$
$$D=(6-\cos\alpha_i-\cos\beta_j) -\sqrt{(6-\cos\alpha_i-\cos\beta_j)^2-4(7-3\cos\alpha_i-3\cos\beta_j-cos\alpha_i\cos\beta_j  )}.$$

Considering that $m,n$ approach infinity, we have
\begin{eqnarray}
\nonumber && \lim_{m\to \infty}\lim_{n\to \infty}\frac{
\mathscr{LEL}\Big(UJL(n,m)\Big)}{2mn} \\
&=&\frac{1}{8\pi^2} \int_{0}^{2\pi}\int_{0}^{2\pi} \sqrt{C} \cdot
\ dx dy +\frac{1}{8\pi^2}\int_{0}^{2\pi}\int_{0}^{2\pi}
\sqrt{D} \cdot \ dx dy \\
   &\approx & 4.6841.
\end{eqnarray}

Consequently, according to the equality (9), we can get the
asymptotic Laplacian-energy-like invariant
$\mathscr{LEL}\Big(UJL(n,m)\Big)\approx 9.3682 mn,$ as $m,n\to
\infty$ .

Summing up, we complete the proof. \qed

\vspace{5pt} \noindent {\bf Remark 5.3}\ The numerical integration
value in equality (8) is calculated with MATLAB software
calculation, i.e.,
$$\frac{1}{8\pi^2}\int_{0}^{2\pi}\int_{0}^{2\pi} \sqrt{C}\cdot\ dx
dy\approx 2.9874,~~~
\frac{1}{8\pi^2}\int_{0}^{2\pi}\int_{0}^{2\pi} \sqrt{D}\cdot\ dx
dy\approx 1.6967. $$

\newpage
 \vspace{5pt} \noindent
{\bf Acknowledgments}\

The work of J. B. Liu is partly supported by the Natural Science
Foundation of Anhui Province of China under Grant No. KJ2013B105
and the National Science Foundation of China under Grant
Nos.11471016 and 11401004. The work of X. F. Pan is partly
supported by the National Science Foundation of China under Grant
Nos. 10901001, 11171097, and 11371028.

\end{document}